\documentclass[11pt]{amsart}
\usepackage{amsmath}
\usepackage{amssymb}
\usepackage{amsthm}
\usepackage{latexsym}
\usepackage{hyperref}
\usepackage{enumerate}

\newtheorem{thm}{Theorem}[section]
\newtheorem{cor}[thm]{Corollary}
\newtheorem{lem}[thm]{Lemma}
\newtheorem{prop}[thm]{Proposition}

\newtheorem{cond}[thm]{Condition}

\providecommand{\norm}[1]{\left\| #1 \right\|}
\newcommand{\enuma}[1]{\begin{enumerate}[\textup{(}a\textup{)}] {#1} \end{enumerate}}
\DeclareMathOperator{\Mod}{Mod}

\newcommand{\mh}{\mathbb}
\newcommand{\mr}{\mathrm}
\newcommand{\mc}{\mathcal}

\newcommand{\scs}{\scriptstyle}

\newcommand{\N}{\mathbb N}
\newcommand{\Z}{\mathbb Z}

\newcommand{\R}{\mathbb R}
\newcommand{\C}{\mathbb C}
\newcommand{\ep}{\epsilon}

\begin{document}

\title{Resolutions of tempered representations of reductive $p$-adic groups}
\author{Eric Opdam}
\address{Korteweg-de Vries Institute for Mathematics\\
Universiteit van Amsterdam\\
Science Park 904\\
1098 XH Amsterdam\\
The Netherlands}
\email{e.m.opdam@uva.nl}
\author{Maarten Solleveld}
\address{Institute for Mathematics, Astrophysiscs and Particle Physiscs\\
Radboud Universiteit Nijmegen\\
Heyendaalseweg 135\\
6525AJ Nijmegen\\
The Netherlands}
\email{m.solleveld@science.ru.nl}
\date{\today}
\subjclass[2010]{22G25, 22E50}

\maketitle

\begin{abstract}
Let $G$ be a reductive group over a non-archimedean local field and let $\mc S (G)$
be its Schwartz algebra. We compare Ext-groups of tempered $G$-representations in several
module categories: smooth $G$-representations, algebraic $\mc S (G)$-modules, bornological
$\mc S (G)$-modules and an exact category of $\mc S (G)$-modules on LF-spaces, 
which contains all admissible $\mc S (G)$-modules.
We simplify the proofs of known comparison theorems for these Ext-groups, due to
Meyer and Schneider--Zink. Our method is based on the Bruhat--Tits building of $G$ and
on analytic properties of the Schneider--Stuhler resolutions.
\end{abstract}

\tableofcontents

\section*{Introduction}
Let $G$ be the group of $\mathbb F$-rational points of a connected reductive linear 
algebraic group defined over a non-archemedean local field $\mathbb F$ of arbitrary characteristic. 
Let $\mc H(G)$ denote the Hecke algebra of locally constant compactly supported complex 
functions on $G$, and $\mc S(G)$ the Harish-Chandra Schwartz algebra of $G$.
The abelian category Mod$(G)$ of complex smooth representations of $G$ is equivalent to 
the category Mod$(\mc H(G))$ of nondegenerate $\mc H(G)$-modules. By \cite[Proposition 1]{ScSt2} 
an admissible representation $V$ of $G$ is tempered if and only if it 
extends to a module of $\mc S(G)$, and then such an extension is unique. Let $V,W$ 
be $\mc S(G)$-modules, with $V$ admissible. A profound theorem due to Schneider and 
Zink \cite{ScZi2} (based on work of R. Meyer) states that for all $n \in \Z_{\geq 0}$:
\begin{equation}\label{eq:intro1}
\textup{Ext}^n_{\mc H(G)}(V,W) = \textup{Ext}^n_{\mc S(G)}(V,W) .
\end{equation}
If $W$ is also admissible, then both $V$ and $W$ admit a canonical structure as 
LF-spaces such that they become complete topological modules over the 
LF-algebra $\mc S(G)$. We introduce an exact category $\textup{Mod}_{LF}(\mc S(G))$ 
of certain LF-modules over $\mc S(G)$, whose exact sequences are split as 
LF-spaces, and then one also has
\begin{equation}\label{eq:intro2}
\textup{Ext}^n_{\mc S(G)}(V,W) = \textup{Ext}^n_{\Mod_{LF}(\mc S(G))}(V,W) .
\end{equation}
One can choose a good compact open subgroup $K$ such that $V=\mc H(G)V^K$ 
and $W=\mc H(G)W^K$. Now $V^K$ and $W^K$ are finite dimensional modules over 
the $K$-spherical Hecke algebra $\mc H(G,K) := e_K\mc H(G) e_K$ which uniquely extend 
to topological modules over the Fr\'echet algebra $\mc S(G,K) := e_K\mc S(G) e_K$.  
In that context we have  
\begin{equation}\label{eq:intro3}
\textup{Ext}^n_{\mc H(G,K)}(V^K,W^K)=
\textup{Ext}^n_{\Mod_{Fr}(\mc S(G,K))}(V^K,W^K) ,
\end{equation}
where $\Mod_{Fr}$ denotes the exact category of Fr\'echet modules
with linearly split exact sequences.

These are powerful statements wich provide a link between harmonic analysis and 
homological properties of admissible smooth representations of $G$. For example it 
follows that a discrete series representation of $G$ is a projective module  
in the full subcategory of $\mc H(G)$-modules which are restrictions of 
$\mc S(G)$-modules. The identities \eqref{eq:intro1}, \eqref{eq:intro2} and \eqref{eq:intro3} 
were used in \cite{OpSo3} to explicitly compute the spaces $\textup{Ext}^i_{\mc H(G)}(V,W)$
for irreducible tempered admissible representations of $G$ in terms of  
analytic R-groups. As a further consequence, we proved the Kazhdan orthogonality 
relations for admissible characters of $G$ directly from the Plancherel isomorphism 
for $\mc S(G)$. These applications motivated us to revisit the proofs of the results of 
Meyer and the subsequent results of Schneider, Stuhler and Zink discussed above.

Equation \eqref{eq:intro1} is somewhat unexpected since $\mc S(G)$ is not a flat ring over 
$\mc H(G)$. Meyer's proof of these type of results \cite{Mey2} relies
in an essential way on the machinery of bornological vector spaces. In the present paper 
we prove the results in a way which is intuitively more clear and which reveals their geometric 
origin. The methods we are using are similar to those used in \cite{OpSo1} for the analogous 
statements for tempered modules over an affine Hecke algebra. The pleasant surprise is that
such an explicit construction of a continuous contraction of the Schneider--Stuhler resolutions 
is still possible in this more complicated context, and the computations are not too unpleasant.
 
First recall the construction of Schneider and Stuhler of a functorial projective 
resolution $C_*(\mc B(G),V)$ of $V$ by $G$-equivariant sheaves on the Bruhat--Tits building 
$\mc B(G)$. We start by constructing a functorial contraction of the $K$-invariant part 
$V^K \leftarrow C_*(\mc B(G),V)^K$ of this resolution of $V$, where $K$ runs over a 
neighborhood basis of $G$ consisting of good compact open subgroups. This is a projective 
resolution of $V^K$ as a $\mc H(G,K)$-module.
The construction of these contractions reflects the contractibility of the affine building 
$\mc B(G)$. Next we show directly that these contractions extend continuously to 
the natural Fr\'echet completion $C_*^t(\mc B(G),V)^K$ of this 
resolution of $V^K$. This shows that the Fr\'echet completion of the  
resolution is an admissible projective resolution of the $\mc S(G,K)$-module $V^K$ in  
$\textup{Mod}_{Fr}(\mc S(G,K))$, and leads to \eqref{eq:intro3}.
 
Given a good maximal compact subgroup $K$, we denote by Mod$(\mc H(G),K)$ the 
full subcategory of $\mc H(G)$-modules $V$ such that $V = \mc H(G)V^K$.
By well known results of Bernstein the functor from Mod$(\mc H(G),K)$ to 
Mod$(\mc H(G,K))$ given by $V\to V^K$ is an equivalence of categories, and by 
results of Schneider and Zink \cite{ScZi1} a similar statement holds
for modules $V$ over $\mc S(G)$ satisfying $V=\mc S(G)V^K$. 
If we take $K$ sufficiently small, such that $V^K$ generates $V$ as a $G$-module  
and $W^K$ generates $W$ as a $G$-module, then one derives \eqref{eq:intro1} 
and \eqref{eq:intro2} from \eqref{eq:intro3} using these equivalences.\\[1mm]

\textbf{Acknowledgements.}
It is a pleasure to thank Joseph Bernstein for stimulating discussions.
During this research Eric Opdam was supported by ERC advanced grant no. 268105.
\vspace{2mm}

\section{The differential complex}

Let $\mh F$ be a non-archimedean local field of arbitrary characteristic. Let $\mc G$ be a connected
reductive algebraic group defined over $\mh F$ and let $G = \mc G (\mh F)$ be its $\mh F$-rational 
points. We briefly call $G$ a reductive $p$-adic group. Let $Z(\mc G)$ be the centre of $\mc G$ and
denote by $X_* (H)$ the set of $\mh F$-algebraic cocharacters of an $\mh F$-group $H$.

The (enlarged) Bruhat--Tits building of $G$ is
\begin{equation}
\mc B (G) = \mc B (\mc G, \mh F) = \mc B (\mc G / Z(\mc G),\mh F) \times 
X_* \big( Z(\mc G (\mh F)) \big) \otimes_\Z \R .
\end{equation}
Recall that $\mc B (\mc G / Z(\mc G),\mh F)$ is in a natural way a polysimplicial complex 
with a $G$-action. The choice of a basis of the lattice $X_* ( Z(G) )$ 
induces a polysimplicial structure on $X_* ( Z(G)) \otimes_\Z \R$, 
isomorphic to a direct product of some copies of $\R$ with the intervals $[n,n+1]$ as 1-simplices. 
The resulting polysimplicial structure on $\mc B (G)$ is $G$-stable because $G$ acts on 
$X_* ( Z(G)) \otimes_\Z \R$ via translation over $X_* ( Z(G))$. 

A crucial role will be played by a system of compact open subgroups of $G$ introduced by
Schneider and Stuhler \cite{ScSt}. The group associated to a given polysimplex $\sigma \subset \mc B (G)$
and a natural number $e$ is denoted $U_\sigma^{(e)}$. We will need the following properties,
which can be found in \cite[Chapter 1]{ScSt} and in \cite[Theorem 5.5]{MeSo2}.

\begin{prop}\label{prop:1.1}
\enuma{
\item $U_\sigma^{(e)}$ is an open pro-$p$ subgroup of the stabilizer $G_\sigma$ of $\sigma$.
\item The collection $\{U_\sigma^{(e)} \mid e \in \N \}$ is a neighborhood basis of $1$ in $G$.
\item $U_\sigma^{(e)}$ depends only on the projection of $\sigma$ on $\mc B (\mc G / Z(\mc G),\mh F)$.
\item $g U_\sigma^{(e)} g^{-1} = U_{g \sigma}^{(e)}$ for all $g \in G$, in particular $U_\sigma^{(e)}$
is a normal in $G_\sigma$.
\item $U_\sigma^{(e)}$ fixes the star of $\sigma$ in $\mc B (G)$ pointwise.
\item $U_\sigma^{(e)}$ is the product (in any order) of the groups $U_x^{(e)}$, where $x$ runs over
the vertices of $\sigma$.
\item If $\sigma_1 ,\sigma_2$ and $\sigma_3$ are polysimplices of $\mc B (G)$ such that $\sigma_2$
lies in every apartment of $\mc B (G)$ that contains $\sigma_1 \cup \sigma_3$, then
$U_{\sigma_2}^{(e)} \subset U_{\sigma_1}^{(e)} U_{\sigma_3}^{(e)}$.
}
\end{prop}

Let $\mc O$ be the ring of integers of $\mh F$ and let $\pi$ be a uniformizer of $\mc O$. Let $p$ be
the characteristic of the residue field $\mc O / \pi \mc O$. We recall the main result of \cite{MeSo1},
which works in the generality of modules over a commutative unital ring $R$ in which $p$ is invertible.
Let $C_c^\infty (G,R)$ be the $R$-module of locally constant, compactly supported functions $G \to R$.
Since $G$ is locally a pro-$p$ group, there exists a Haar measure on $G$ such that all pro-$p$ 
subgroups of $G$ have volume in $p^\Z$. We fix such a measure once and for all. Thus we obtain a 
convolution product on $C_c^\infty (G;R)$, which makes it into an $R$-algebra denoted $\mc H (G;R)$. 
Let Mod$(\mc H (G;R))$ be the category of $\mc H (G;R)$-modules $V$ with $\mc H (G;R) V = V$. It is 
naturally equivalent to the category $\Mod_R (G)$ of smooth $G$-representations on $R$-modules.

Now we describe how the above objects can be used to construct resolutions of certain modules. Given
any polysimplex $\sigma$, let $e_\sigma = e_{U_\sigma^{(e)}}$ be the corresponding idempotent of
$\mc H (G;R)$; it exists because the volume of $U_\sigma^{(e)}$ is invertible in $R$. For any
$V \in \Mod_R (G) ,\; e_\sigma V = V^{U_\sigma^{(e)}}$ is the $R$-submodule of 
$U_\sigma^{(e)}$-invariant elements. For any polysimplicial subcomplex $\Sigma \subset \mc B (G)$
let $\Sigma^{(n)}$ be the collection of $n$-dimensional polysimplices of $\Sigma$. We put
\[
C_n (\Sigma ; V) := \bigoplus_{\sigma \in \Sigma^{(n)}} R \sigma \otimes_R e_\sigma V .
\]
We fix an orientation of the polysimplices of $\mc B (G)$ and we identify $-\sigma$ with $\sigma$
oriented in the opposite way. This allows us to write the boundary of $\sigma$ in the polysimplicial 
sense \cite[Section 2.1]{OpSo1} as 
\[
\partial \sigma = \sum\nolimits_\tau [\sigma : \tau] \tau \quad 
\text{with} \quad [\sigma : \tau] \in \{1,0,-1\} .
\]
We have $[\sigma : \tau] = 0$ unless $\tau \subset \sigma$, and in that case 
Proposition \ref{prop:1.1}.f tells us that $U_\tau^{(e)} \subset U_\sigma^{(e)}$ and 
$e_\tau V \supset e_\sigma V$. Thus we can define a differential
\begin{align*}
& \partial_n : C_n (\Sigma ; V) \to C_{n-1} (\Sigma ; V) , \\
& \partial_n (\sigma \otimes v) = \partial \sigma \otimes v = 
\sum\nolimits_\tau [\sigma :\tau] \tau \otimes v
\end{align*}
and an augmentation
\begin{align*}
& \partial_0 : C_0 (\Sigma ; V) \to V = C_{-1} (\Sigma ; V) ,\\
& \partial_0 (x \otimes v) = v .
\end{align*}
Since $\partial^2 = 0 ,\; \big( C_* (\Sigma;V),\partial_* \big)$ is a differential complex. 
The group $G$ acts on $C_n (\mc B (G);V)$ by
\begin{equation}
g (\sigma \otimes v) = g \sigma \otimes g v ,
\end{equation}
where $g \sigma$ is endowed with the orientation that makes $g : \sigma \to g \sigma$ orientation 
preserving. Clearly $\partial_*$ is $G$-equivariant, so $\big( C_* (\mc B (G);V) , \partial_* \big)$
is a complex of $\mc H (G;R)$-modules. 

\begin{thm}\label{thm:1.2}
Let $\Sigma \subset \mc B (G)$ be convex.
\enuma{
\item The differential complex $\big( C_* (\Sigma;V),\partial_* \big)$ is acyclic and $\partial_0$
induces a bijection
\[
H_0 \big( C_* (\Sigma;V),\partial_* \big) \to \sum_{x \in \Sigma^{(0)}} V^{U_x^{(e)}} .
\]
\item If $V = \sum_{x \in \mc B (G)^{(0)}} V^{U_x^{(e)}}$, then $\big( C_* (\mc B (G);V) , 
\partial_* \big)$ is a resolution of $V$ in $\Mod_R (G)$. This resolution is projective if
the order of $G_\sigma / U_\sigma^{(e)}$ is invertible in $R$ for every polysimplex $\sigma$.
}
\end{thm}
\emph{Proof.}
(a) is \cite[Theorem 2.4]{MeSo1}. Although in \cite{MeSo1} the affine building of $\mc G / Z (\mc G)$
is used, this does not make any difference for the proof. In particular the crucial
\cite[Theorem 2.12]{MeSo1} is also valid in our setup. \\
(b) The special case where $R = \C$ and $Z(G)$ is compact was proven in \cite[Theorem II.3.1]{ScSt}.
It remains to show that $C_n (\mc B (G);V)$ is projective under the indicated conditions. 
Let $\sigma_1, \ldots, \sigma_d$ be representatives for the $G$-orbits of $n$-dimensional 
polysimplices in $\mc B (G)$ and let $\ep_{\sigma_i} : G_{\sigma_i} \to \{ 1,-1\}$ be the orientation
character of $\sigma_i$. By construction
\begin{equation}\label{eq:1.2}
C_n (\mc B (G);V) = 
\bigoplus_{i=1}^d \mr{ind}_{G_{\sigma_i}}^G (\ep_{\sigma_i} \otimes e_{\sigma_i} V) = 
\bigoplus_{i=1}^d \mr{ind}_{R [G_{\sigma_i} / U_{\sigma_i}^{(e)}]}^{\mc H (G;R)} 
\big (\ep_{\sigma_i} \otimes V^{U_{\sigma_i}^{(e)}} \big) .
\end{equation}
By assumption $|G_{\sigma_i} / U_{\sigma_i}^{(e)} | \in R^\times$, so the category of 
$R [G_{\sigma_i} / U_{\sigma_i}^{(e)}]$-modules is semisimple. In particular 
$\ep_{\sigma_i} \otimes V^{U_{\sigma_i}^{(e)}}$ is projective in this category, which by Frobenius
reciprocity implies that $\mr{ind}_{R [G_{\sigma_i} / U_{\sigma_i}^{(e)}]}^{\mc H (G;R)} 
\big( \ep_{\sigma_i} \otimes V^{U_{\sigma_i}^{(e)}} \big)$ is projective in 
Mod$(\mc H (G;R)). \qquad \Box$ 
\vspace{2mm}

\section{A functorial contraction}

From now on we fix $e \in \N$ and a special vertex $x_0 \in \mc B (G)$. By Proposition 
\ref{prop:1.1}.d $K = U_{x_0}^{(e)}$ is a normal subgroup of the good maximal compact subgroup 
$G_{x_0}$ of $G$. Let $\mc H (G,K;R)$ be the subalgebra of $\mc H (G;R)$ consisting of all 
$K$-biinvariant elements and let $\Mod_R (G,K)$ be the full subcategory of $\Mod_R (G)$ 
made of all objects $V$ for which $V = \mc H(G;R) V^K$.

\begin{thm}\label{thm:2.1}
Suppose that $R$ is an algebraically closed field whose characteristic is banal for $G$, that is,
does not divide the pro-order of any compact subgroup of $G$. 
\enuma{
\item The exact functor 
\[
\begin{array}{ccc}
\Mod_R (G) & \to & \Mod(\mc H (G,K;R)) , \\
V & \mapsto & e_K V = V^K 
\end{array}
\]
provides an equivalence of categories $\Mod_R (G,K) \to \Mod(\mc H (G,K;R))$, 
with quasi-inverse $W \mapsto \mc H (G;R) e_K \otimes_{\mc H (G,K;R)} W$.
\item $\Mod_R (G,K)$ is a direct factor of the category $\Mod_R (G)$, 
so the above functors preserve projectivity.
}
\end{thm}
\emph{Proof.}
For $R = \C$ this is due to Bernstein, see \cite{BeDe}. Vign\'eras \cite{Vig2} observed that 
Bernstein's proof remains valid for $R$ as indicated. $\qquad \Box$ \\[2mm]

We remark that it is likely that Theorem \ref{thm:2.1} is valid for much more general rings $R$.
By \cite[Section 3]{MeSo1} this is the case for somewhat different idempotents of $\mc H (G;R)$.

\begin{thm}\label{thm:2.2}
\enuma{
\item The augmented complex
\[
V^K \leftarrow C_* (\mc B (G);V)^K
\]
admits a contraction which is natural in $V \in \Mod_R (G,K)$.
\item Suppose that $R$ is an algebraically closed field whose characteristic is banal for $G$. Then
\[
\mc H (G,K;R) \leftarrow C_* (\mc B (G);\mc H (G;R) e_K )^K
\]
is a projective $\mc H(G,K;R)$-bimodule resolution which admits a right $\mc H(G,K;R)$-linear 
contraction.
}
\end{thm}
\emph{Remark.} $C_* (\mc B (G);V)$ need not be generated its $K$-invariant vectors. The point is
that there can be polysimplices $\sigma$ that do not contain any element of $G x_0$, and then
$U_\sigma^{(e)}$ does not contain any conjugate of $K$.

\emph{Proof.}
(a) In the lowest degree it is easy, for $v \in V^K = C_{-1}(\mc B (G);V)^K$ we put
\[
\gamma_{-1}(v) = x_0 \otimes v.
\]
In higher degrees, consider an apartment $A$ of $\mc B (G)$ containing $x_0$ and let $T$ be the
associated maximal split torus of $G$. In \cite[Section 2.1]{OpSo1} a contraction $\gamma$ of the
augmented differential complex $C_* (A;\Z)$ is constructed, with the properties:
\begin{enumerate}
\item $\gamma$ is equivariant for the action of the Weyl group $N_{G_{x_0}}(T) / (T \cap G_{x_0})$
on $A$.
\item For any polysimplex $\sigma \subset A$, the support of $\gamma (\sigma)$ is contained 
in the hull of $\sigma \cup \{x_0\}$. 
\end{enumerate}
Recall that the hull of a subset $X$ of a thick affine building (as in property 1) is the 
intersection of all apartments that contain $X$. It is the best polysimplicial approximation
of the convex closure of $X$. 
An alternative, more explicit, description is given in \cite[Section 1.1]{OpSo1}. This works well for
the Bruhat--Tits building of $G / Z(G)$, which is the setting of \cite{OpSo1}.

For $X_* (Z(G)) \otimes_\Z \R$ an intersection of apartments is not suitable because there is only
one apartment. But by means of a basis of $X_* (Z(G))$ we have already 
identified this lattice with $\Z^d$. We define the hull of any subset $X \subset A$ as the hull of the 
projection of $X$ on $\mc B (G/Z(G))$ times the smallest box 
\[
\prod\nolimits_{i=1}^d [n_i,n'_i] \subset X_* (Z(G)) \otimes_\Z \R \quad 
\text{with} \quad n_i,n'_i \in \Z
\]
that contains the projection of $X$ on $X_* (Z(G)) \otimes_\Z \R$. To reconcile this with 
\cite{OpSo1}, let $\{\beta_1, \ldots, \beta_d\}$ be the dual basis of $X^* (Z(G))$
and regard $\{\beta_1 ,-\beta_1 ,\ldots, \beta_d ,-\beta_d\}$ as a root system of 
type $A_1^d$ in $X^* (Z(G)) \otimes_\Z \R$. Then our hull agrees with the description 
given in \cite[Section 1.1]{OpSo1}, so the construction of $\gamma$ applies to $A$.

For an elementary tensor $\sigma \otimes v_\sigma$ we define
\[
\gamma_n (\sigma \otimes v_\sigma) = 
\gamma (\sigma) \otimes v_\sigma \in C_{n+1}(A;\Z) \otimes_\Z V.
\]
By property 2 this does not depend on the choice of the apartment $A$ and it clearly is functorial
in $V$. Recall from \cite[Proposition 7.4.8]{BrTi1} that, whenever $B \subset A, g \in G$ and
$g(B) \subset A$, there exists $n \in N_G (T)$ such that $n(b) = g(b)$ for all $b \in B$. 
With property 1 it follows that $\gamma_n$ extends to a $G_{x_0}$-equivariant map
\[
\gamma_n : C_n (\mc B (G);V) \to C_{n+1}(\mc B (G);\Z) \otimes_\Z V.
\]
Consider a typical $K$-invariant element $e_K (\sigma \otimes v_\sigma) \in C_n (\mc B (G);V)^K$.
By the $K$-equivariance of $\gamma_n$
\begin{equation}\label{eq:2.6}
\gamma_n (e_K (\sigma \otimes v_\sigma )) = 
e_K (\gamma (\sigma) \otimes v_\sigma ) = e_K \big( \sum\nolimits_\tau \gamma_{\sigma \tau} 
\tau \otimes v_\sigma \big) \in C_{n+1} (\mc B (G);V)^K
\end{equation}
By property 2 the support of \eqref{eq:2.6} lies in the $K$-orbit of the hull of $\{x_0\} \cup \sigma$.
We fix a polysimplex $\tau_0$ in this support and we have a closer look at the 
value of \eqref{eq:2.6} at $\tau_0$. Let us write $G_\tau^+$ for the subgroup of $G_\tau$ consisting
of the elements that preserve the orientation of the polysimplex $\tau$. Then the value of
\eqref{eq:2.6} at $\tau_0$ is
\begin{equation}\label{eq:2.2}
\sum_{\tau \subset A \cap K \tau_0} \sum_{k \in K / (K \cap G_\tau^+), k \tau = \pm \tau_0}
\gamma_{\sigma\tau} \epsilon (k,\tau) k \frac{e_{K \cap G_\tau^+}}{[K : K \cap G_\tau^+]} v_\sigma ,
\end{equation}
where $\epsilon (k,\tau) \in \{\pm 1\}$ is defined by $k\tau = \epsilon (k,\tau) \tau_0$. 
By Proposition \ref{prop:1.1}.g we have $U_\tau^{(e)} \subset U_{x_0}^{(e)} U_\sigma^{(e)}$
for all $\tau$ occurring in the above sum. 
However, we need a more precise version. Let $\Phi$ be the root system of $(G,T)$ and let 
$\Phi^+$ be a system of positive roots such that $\sigma$ lies in the positive Weyl chamber $A^+$. 
Let $U^+$ and $U^-$ be the unipotent subgroups of $G$ associated to $\Phi^+$ and $-\Phi^+$. 
The constructions in \cite[Section 1.2]{ScSt} entail that 
\begin{equation}\label{eq:2.7}
\begin{array}{lclcl}
U_{x_0}^{(e)} \cap U^+ & \subset & U_\tau^{(e)} \cap U^+ & \subset & U_\sigma^{(e)} \cap U^+, \\
U_{x_0}^{(e)} \cap Z_G (T) & = & U_\tau^{(e)} \cap Z_G (T) & = & U_\sigma^{(e)} \cap Z_G (T), \\
U_{x_0}^{(e)} \cap U^- & \supset & U_\tau^{(e)} \cap U^- & \supset & U_\sigma^{(e)} \cap U^- .
\end{array}
\end{equation}
By \cite[Corollary I.2.8]{ScSt} 
\begin{equation}\label{eq:2.4}
U_\tau^{(e)} = (U_\tau^{(e)} \cap U^-) (U_\tau^{(e)} \cap Z_G (T)) (U_\tau^{(e)} \cap U^+) ,
\end{equation}
which is contained in
\[
\begin{split}
& (U_\tau^{(e)} \cap U^-) (U_\sigma^{(e)} \cap Z_G (T)) (U_\sigma^{(e)} \cap U^+) \; \subset \\
& (U^- \cap G_\tau^+ \cap U_{x_0}^{(e)}) (U_\sigma^{(e)} \cap Z_G (T) U^+) 
\; \subset \; (G_\tau^+ \cap K) U_\sigma^{(e)} .
\end{split}
\]
By assumption $v_\sigma \in V^{U_\sigma^{(e)}}$, so the above means that 
\[
ke_{K \cap G_\tau^+} v_\sigma = k e_{K \cap G_\tau^+} e_\tau v_\sigma = 
k e_\tau e_{K \cap G_\tau^+} v_\sigma = 
e_{\tau_0} k e_{K \cap G_\tau^+} v_\sigma \in  V^{U_{\tau_0}^{(e)}} ,
\]
where we used that $G_\tau$ normalizes $U_\tau^{(e)}$ for the second equality.
Consequently \eqref{eq:2.2} lies in $V^{U_{\tau_0}^{(e)}}$ and
\begin{equation}\label{eq:2.13}
\gamma_n (e_K (\sigma \otimes v_\sigma )) \in C_{n+1}(\mc B (G);V)^K .
\end{equation}
Because $\gamma \partial + \partial \gamma =$ id on $C_* (A;\Z) ,\; \gamma_*$ is a contraction
of the augmented complex $\big( C_* (\mc B (G),V), \partial_* \big)$. \\
(b) It only remains to show that $C_n (\mc B (G);\mc H (G;R) e_K )^K$ is projective as a
$\mc H (G,K;R)$-bimodule. By \eqref{eq:1.2} it equals
\begin{equation}\label{eq:2.1}
\bigoplus_{i=1}^d \big( \mr{ind}_{G_{\sigma_i}}^G (\ep_{\sigma_i} \otimes e_{\sigma_i} 
\mc H (G;R) e_K )\big)^K = 
\bigoplus_{i=1}^d e_K \mc H (G;R) \underset{\scs{\mc H (G_{\sigma_i};R)}}{\otimes} \ep_{\sigma_i} 
\otimes e_{\sigma_i} \mc H (G;R) e_K .
\end{equation}
Consider the $\mc H(G;R)$-bimodule $\mc H (G;R) \otimes_{\mc H (G_\sigma ;R)} \ep_\sigma \otimes 
e_\sigma \mc H (G;R)$, for any polysimplex $\sigma$. Since $U_\sigma^{(e)}$ is normalized by $G_\sigma$ 
and fixes $\sigma$ pointwise, $e_\sigma\in\mc H(G_\sigma;R)$ is a central idempotent and
$e_\sigma\epsilon_\sigma=\epsilon_\sigma$. As $R$ has banal characteristic, the group algebra 
\begin{equation}\label{eq:2.5}
R [G_\sigma / U_\sigma^{(e)}] = \mc H(G_\sigma,U_\sigma^{(e)};R) = e_\sigma\mc H(G_\sigma;R) 
\end{equation}
is a semisimple direct summand of $\mc H(G_\sigma;R)$. Given an irreducible 
representation $\rho$ of \eqref{eq:2.5} let $e_\rho$ be the corresponding central idempotent 
and $p_\rho \in e_\rho R [G_\sigma / U_\sigma^{(e)}]$ an idempotent of rank 1. 
Then $\sum_\rho (p_{\ep_\sigma\otimes\rho} \otimes p_\rho)$ is an 
idempotent in $\mc H (G;R) \otimes_R \mc H (G;R)^{op}$ and as $\mc H (G;R)$-bimodules
\begin{align}
& \nonumber \mc H (G;R) \underset{\scs{\mc H (G_\sigma ;R)}}{\otimes} \ep_\sigma 
\otimes e_\sigma \mc H (G;R) = \\
& \label{eq:2.3} \bigoplus\nolimits_\rho \mc H (G;R) e_{\epsilon_\sigma\otimes\rho} 
\underset{\scs{\mc H (G_\sigma ;R)}}{\otimes} e_\rho \mc H(G;R) = \\
& \nonumber \bigoplus\nolimits_\rho \mc H (G;R) p_{\epsilon_\sigma\otimes\rho}  
\otimes_R p_\rho \mc H(G;R) = \\
& \nonumber \big( \mc H (G;R) \otimes_R \mc H (G;R)^{op} \big) 
\sum\nolimits_\rho (p_{\epsilon_\sigma\otimes\rho}  \otimes p_\rho) .
\end{align}
In particular $\mc H (G;R) \otimes_{\mc H (G_\sigma ;R)} \ep_\sigma \otimes e_\sigma \mc H (G;R)$
is projective as a $\mc H(G;R)$-bimodule. By Theorem \ref{thm:2.1}, applied to the category of
smooth $G \times G^{op}$-representations on $R$-modules, 
\[
e_K \mc H (G;R) \underset{\scs{\mc H (G_\sigma ;R)}}{\otimes} 
\ep_\sigma \otimes e_\sigma \mc H (G;R) e_K
\]
is projective in Mod$(\mc H (G,K;R) \otimes_R \mc H (G,K;R)^{op})$. Together with \eqref{eq:2.1}
this completes the proof. $\qquad \Box$ 
\\[2mm]

We record some useful properties of the contraction constructed above.
\begin{cor}\label{cor:2.3}
The contraction $\gamma$ from Theorem \ref{thm:2.2}.a satisfies:
\enuma{
\item $\gamma$ is $G_{x_0}$-equivariant.
\item For any polysimplex $\sigma$, the support of $\gamma (\sigma)$ is contained 
in the hull of $\sigma \cup \{x_0\}$. In particular it intersects every $G_{x_0}$-orbit in at
most one polysimplex.
\item The exists $M_\gamma \in \N$ such that $\gamma (\sigma) = \sum_\tau \gamma_{\sigma \tau} 
\tau$ with $|\gamma_{\sigma \tau}| \leq M_\gamma$ for all polysimplices $\sigma$ and $\tau$ of 
$\mc B (G)$.
}
\end{cor}
\emph{Proof.}
(a) and the first half of (b) follow from the properties of the contraction of $C_* (A;\Z)$ that we
used in the proof of Theorem \ref{thm:2.2}. Let $A'$ be an apartment containing $\sigma$ and
$x_0$, corresponding to a maximal split torus $T'$. Then the hull of $\sigma \cup \{x_0\}$ is
contained in a Weyl chamber for the Weyl group $N_{G_{x_0}}(T') / (T' \cap G_{x_0})$. Two points
of $A'$ are in the same $G_{x_0}$-orbit if and only if they are in one
$N_{G_{x_0}}(T') / (T' \cap G_{x_0})$-orbit. Hence this hull, and in particular the support of 
$\gamma (\sigma)$, intersects every $G_{x_0}$-orbit in at most one polysimplex.\\
(c) is a direct consequence of the definition of $\gamma$ and of the corresponding
property of the contraction from \cite[Section 2.1]{OpSo1}. $\qquad \Box$
\\[2mm]

Let $d$ be a $G$-invariant metric on $\mc B (G)$. Then the restriction of $d$ to an apartment
$A$ comes from an inner product on $A$. We may and will assume that:
\begin{itemize}
\item in any apartment $A ,\; X_* (Z(G)) \otimes_\Z \R$ is orthogonal to $A \cap \mc B (G / Z(G))$;
\item the chosen basis of $X_* (Z(G)) \otimes_\Z \R$ is orthogonal.
\end{itemize}
The following elementary property of the hull of a polysimplex of 
$\mc B(G)$ will be used in the proof of Theorem \ref{thm:4.2}.

\begin{lem}\label{lem:2.8}
There exists $\delta \in \mc B (G)$ such that, for every polysimplex $\sigma$, the diameter
of the hull of $\sigma \cup \{x_0\}$ is at most $d(x_0,\sigma) + d(x_0,\delta)$. 
\end{lem}
\emph{Proof.}
Denote the hull of $\sigma \cup \{x_0\}$ by $\mc H$. Let the apartment $A$ and the 
root system $R = R(G,T) \cup \pm \{\beta_1 ,\ldots, \beta_d\}$ be as in the proof of 
Theorem \ref{thm:2.2}. The above assumptions on $d$ imply that there exist $c_\alpha > 0$
such that
\begin{equation}\label{eq:2.9}
d(x,x_0)^2 = \sum\nolimits_{\alpha \in R} c_\alpha \alpha (x)^2 
\quad \text{for all } x \in A .
\end{equation}
Choose a system of positive roots $R^+$ such that $\sigma$ (and hence $\mc H$) lies 
in the positive Weyl chamber, with $x_0$ as origin. Then
\[
0 \leq \alpha (x) \leq \min_{y \in \sigma} \alpha (y) + 1
\quad \text{for all } x \in \mc H, \alpha \in R^+,
\]
see \cite[\S 1.1]{OpSo1}. Let $\alpha^\vee \in A$ be the coroot of $\alpha$ and put
$\delta = \sum_{\alpha \in R^+} \alpha^\vee / 2$. It is well-known that $\beta (\delta)$
equals the height of $\beta \in R$, which by definition is at least 1 for all $\beta \in R^+$. 
With \eqref{eq:2.9} it follows that, for all $x_1 ,x_2 \in \mc H$:
\begin{align*}
& 0 \leq |\alpha (x_1) - \alpha (x_2) | \leq \min_{y \in \sigma} \alpha (y) + \alpha (\delta)
\quad \text{for all } \alpha \in R^+ , \\
& d(x_1,x_2) \leq d(x_0,\sigma + \delta) \leq d(x_0,\sigma) + d(x_0,\delta) .
\end{align*}
Thus $\delta$ works for this particular $\sigma$. For other polysimplices $\sigma'$ the above 
argument would produce a possibly different $\delta'$. But since all Weyl chambers in all 
apartments containing $x_0$ are conjugate under $G_{x_0} ,\; d(x_0,\delta') = d(x_0,\delta).
\qquad \Box$ 
\vspace{2mm}

\section{Fr\'echet $\mc S (G,K)$-modules}

In this section we consider only complex $G$-representations, so we fix $R = \C$ and we suppress
it from the notation. Recall that the Harish-Chandra--Schwartz algebra $\mc S (G)$ is the union,
over all compact open subgroups $U \subset G$, of the subalgebras 
\[
\mc S (G,U) = e_U \mc S (G) e_U .
\]
By \cite{Vig1} $\mc S (G)$ is the convolution 
algebra of functions $G \to \C$ that are rapidly decreasing with respect to the length 
function $g \mapsto d( g x_0, x_0)$. The subalgebra $\mc S (G,U)$ is nuclear and Fr\'echet. 
Nevertheless $\mc S (G)$ is not a Fr\'echet algebra, it is an inductive
limit of Fr\'echet spaces and its multiplication is only separately continuous. 

Before we move on to Fr\'echet modules, we recall some facts about algebraic $\mc S(G)$-modules.
Let Mod$(\mc S (G))$ denote the category of smooth $\mc S(G)$-modules. Let
Mod$(\mc S(G),U)$ be the subcategory of all $V$ with $V = \mc S (G) V^U$. We write 
$K = U_{x_0}^{(e)}$ as before. The analogue of Theorem \ref{thm:2.1} reads:

\begin{prop}\label{prop:3.1}
\textup{\cite[Lemmas 2.2 and 2.3]{ScZi1}}
\enuma{
\item The exact functor
\[ 
\begin{array}{ccc}
\Mod (\mc S(G)) & \to & \Mod(\mc S (G,K)) , \\
V & \mapsto & e_K V = V^K 
\end{array}
\]
provides an equivalence of categories $\Mod (\mc S (G),K) \to \Mod(\mc S (G,K))$, 
with quasi-inverse $W \mapsto \mc S (G) e_K \otimes_{\mc S (G,K)} W$.
\item Mod$(\mc S(G),K)$ is a direct factor of the category Mod$(\mc S (G))$ and all it objects
$V$ satisfy $V = \mc H (G) V$.
}
\end{prop}

By part (b) the module $\mc S(G) e_K \in \Mod (\mc S (G),K)$ lies also in Mod$(G,K)$. Hence
\begin{equation}\label{eq:3.1}
\mc S(G) e_K = \mc H (G) (\mc S (G) e_K )^K = \mc H (G) e_K \mc S(G) e_K = 
\mc H (G) e_K \underset{\mc H (G,K)}{\otimes} \mc S (G,K) .
\end{equation}
Let $\Mod_{Fr}(\mc S (G,K))$ be the category of Fr\'echet $\mc S (G,K)$-modules, with 
continuous module maps as morphisms. We decree that an exact sequence in this category 
is admissible if it is split exact as a sequence of Fr\'echet spaces. 

Let $F$ be any Fr\'echet space and write the completed projective tensor product as 
$\widehat \otimes$. With respect to the indicated exact structure, modules of the form 
$\mc S (G,K) \widehat \otimes F$ are projective in $\Mod_{Fr}(\mc S (G,K))$. Hence this 
exact category has enough projective objects and all derived functors are well-defined.

Now we specify a category of $\mc S (G)$-modules that is equivalent to 
$\Mod_{Fr}(\mc S (G,K))$. As objects we take all modules $V\in\Mod(\mc S(G),K)$ such that 
$W:=V^K$ belongs to $\Mod_{Fr}(\mc S(G,K))$. Because of Proposition \ref{prop:3.1} we have 
\begin{equation}\label{eq:3.10}
V \cong \mc S (G) e_K \underset{\mc S (G,K)}{\otimes} W \quad
\text{with} \quad W \in \Mod_{Fr}(\mc S (G,K)) . 
\end{equation}
Since $\mc S (G) e_K$ is the union of the subalgebras $e_\sigma^{(r)} \mc S (G) e_K$,
we can write
\begin{equation}\label{eq:3.15}
V = \bigcup_{r \in \N} V^{U_\sigma^{(r)}} \cong 
\bigcup_{r \in \N} e_\sigma^{(r)} \mc S (G) e_K \underset{\mc S (G,K)}{\otimes} W .
\end{equation}
Moreover $e_\sigma^{(r)} \mc S (G) e_K$ is of finite rank as a right $\mc S (G,K)$-module,
so every term $e_\sigma^{(r)} \mc S (G) e_K \underset{\mc S (G,K)}{\otimes} W$ becomes
in a natural way a Fr\'echet space. We endow $V$ with the inductive limit topology
coming from these subspaces, thus making it into an LF-space.

Clearly every $\mc S (G)$-module map $\phi : V \to \tilde V$ sends $V^U$ to $\tilde V^U$,
and by definition $\phi$ is continuous if and only if $\phi |_{V^U}$ is continuous
for every compact open subgroup $U \subset G$. On the other hand, as soon as 
$\phi |_{V^K} : V^K \to \tilde V^K$ is continuous, the aforementioned finite rank
property assures that $\phi |_{V^U}$ is continuous for all $U$. We define 
$\Mod_{Fr}(\mc S (G),K)$ to be the category of all $\mc S (G)$-modules of the
form \eqref{eq:3.10}, with continuous $\mc S(G)$-module maps as morphisms. Exact sequences
in this category are required to be split as sequences of LF-spaces.

Our arguments will also apply to certain modules that are not generated by their 
$U$-invariants for any compact open subgroup $U \subset G$. Let $\Mod_{LF}(\mc S (G))$ 
be the category of all topological $\mc S (G)$-modules $V$ such that:
\begin{itemize}
\item $V^U$ is a Fr\'echet $\mc S (G,U)$-module for every compact open subgroup $U \subset G$;
\item $V$ has the inductive limit topology from the subspaces $V^U$.
\end{itemize}
Of course the morphisms are continuous module maps. As before, we require that exact 
sequences in this category are split as sequences of topological vector spaces. 
In view of \cite[Proposition 1]{ScSt2}, $\Mod_{LF}(\mc S (G))$ naturally contains 
all admissible $\mc S (G)$-modules.
On the other hand, not every LF-space which is a topological $\mc S (G)$-module belongs to 
$\Mod_{LF}(\mc S (G))$, as we require its objects to be LF-spaces in a specific way.
For example, the regular representation on $\mc S (G)$ itself is not in $\Mod_{LF}(\mc S (G))$.

For $V \in \Mod_{Fr}(\mc S (G),K)$ the vector space 
\[
V^U = e_U \mc S (G) e_K \underset{\mc S (G,K)}{\otimes} V^K
\]
has a natural Fr\'echet topology, because $e_U \mc S (G) e_K$ is of finite rank over
$\mc S (G,K)$. Hence $\Mod_{Fr}(\mc S (G),K)$ is contained in $\Mod_{LF}(\mc S (G))$.
In view of Proposition \ref{prop:3.1} and the above considerations we have:

\begin{cor}\label{cor:3.6}
The functors
\[
\begin{array}{ccc}
\Mod_{Fr}(\mc S (G),K) & \longleftrightarrow & \Mod_{Fr}(\mc S (G,K)) \\
V & \mapsto & V^K \\
\mc S (G) e_K \underset{\mc S (G,K)}{\otimes} W & \leftarrow & W
\end{array}
\]
are equivalences of exact categories. Moreover $\Mod_{Fr}(\mc S (G),K)$ is a direct
factor of $\Mod_{LF}(\mc S (G))$.
\end{cor}

For $V \in \Mod_{Fr}(\mc S (G),K)$ \eqref{eq:1.2} and \eqref{eq:3.1} show that
\begin{multline}\label{eq:3.16}
\mc S (G,K) \underset{\mc H (G,K)}{\otimes} C_n (\mc B (G);V)^K =
\bigoplus_{i=1}^d \mc S (G,K) \underset{\mc H (G,K)}{\otimes} e_K \mc H (G) 
\underset{\mc H (G_{\sigma_i})}{\otimes} \ep_{\sigma_i} \otimes e_{\sigma_i} V \\
= \bigoplus_{i=1}^d e_K \mc S (G) \underset{\mc H (G_{\sigma_i})}{\otimes} \ep_{\sigma_i} 
\otimes e_{\sigma_i} V = \bigoplus_{i=1}^d e_K \mc S (G) e_{\sigma_i} 
\underset{\mc H (G_{\sigma_i}; U_{\sigma_i}^{(e)})}{\otimes} \ep_{\sigma_i} \otimes e_{\sigma_i} V .
\end{multline}
By Frobenius reciprocity and Theorem \ref{thm:1.2} this is a projective object of
$\Mod(\mc S (G,K))$. As a vector space it is complete with respect to the projective tensor 
product topology if and only if $e_{\sigma_i} V$ has finite dimension for all $i$. Since 
$V = \mc H (G) V^K$, this happens if and only if $V$ is admissible. 

Unfortunately, for inadmissible $V$ the modules \eqref{eq:3.16} need not form a resolution 
of $V$. The simplest counterexample occurs when $G = \mh F^\times$, a onedimensional torus. 
Then $\mc B (G) \cong \R$ with $\mh F^\times$ acting by translations. Furthermore $K = 1 + 
\omega^{e+1} \mc O$, where $\omega$ is a uniformizer of the ring of integers $\mc O$ of $\mh F$. 
For $V = \mc S (G) e_K \cong \mc S (\Z \times \mc O^\times / K)$ the modules \eqref{eq:3.16} 
form an augmented differential complex
\begin{multline*}
\mc S (\Z \times \mc O^\times / K) \leftarrow \mc S (\Z \times \mc O^\times / K) 
\underset{\mc H (\mc O^\times)}{\otimes} \mc S (\Z \times \mc O^\times / K) \\
\leftarrow \mc S (\Z \times \mc O^\times / K) \underset{\mc H (\mc O^\times)}{\otimes} 
\mc S (\Z \times \mc O^\times / K) .
\end{multline*}
With the appropriate identifications
\[
\partial_0 (f \otimes f') = f f' \quad \text{and} \quad  
\partial_1 (f \otimes f') = f \otimes f' - \omega f \otimes \omega^{-1} f' . 
\]
For any $s \in \mc S (\Z \times \mc O^\times / K) \setminus C_c (\Z \times \mc O^\times / K)$ 
we have $s \otimes e_K - e_K \otimes s \in \ker \partial_0 \setminus \text{im} \, \partial_1$.  
The problem is that, in order to move $s \otimes e_K$ to $e_K \otimes s$, one would need infinitely
long shifts in the direction of $(\omega,\omega^{-1})$, whereas the image of $\partial_1$
can only take care of finitely many such shifts.

To avoid this problem we have to complete \eqref{eq:3.16}, and the best way is to introduce,
for $V \in \Mod_{LF}(\mc S (G))$:
\begin{equation} \label{eq:3.2}
C_n^t (\mc B (G);V) := \bigcup_{r \in \N} \bigoplus_{i=1}^d \big( e_{U_{x_0}^{(r)}} \mc S (G) 
e_{\sigma_i} \underset{\mc H (G_{\sigma_i} ; U_{\sigma_i}^{(e)})}{\widehat \otimes} 
\ep_{\sigma_i} \otimes e_{\sigma_i} V \big) .
\end{equation}
This lies in $\Mod_{LF}(\mc S(G))$ and is a topological completion of $C_n (\mc B (G);V)$. 
Notice the specific order of the operations, which is necessary because (completed) projective 
tensor products do not always commute with inductive limits. We have
\[
C_n^t (\mc B (G);V) = \mc S (G) \otimes_{\mc H (G)} C_n (\mc B (G);V)
\] 
if and only if $V$ is admissible. 

For the augmentation we simply put $C_{-1}^t (V) = V$. 

\begin{lem}\label{lem:3.5}
Let $V \in \Mod_{LF}(\mc S (G))$. The boundary maps $\partial_n$ of $C_* (\mc B (G);V)$ 
extend to continuous $\mc S (G)$-linear boundary maps $\partial_n^t$ of $C_*^t (\mc B (G);V)$.
\end{lem}
\emph{Remark.} An analogous result was used implicitly in \cite[Section 2.3]{OpSo1}. The proof
given here also applies in the setting of \cite{OpSo1}. 

\emph{Proof.} Since the map of Fr\'echet spaces
\[
\begin{array}{ccc}
e_{U_{x_0}^{(r)}} \mc S (G) e_{\sigma_i} \times e_{\sigma_i} V & \to & e_{U_{x_0}^{(r)}} V , \\
(f,v) & \mapsto & f v
\end{array}
\]
is continuous and $\mc H (G_{\sigma_i} / U_{\sigma_i}^{(e)})$-balanced, $\partial_0$ extends
to a continuous map
\begin{equation}\label{eq:3.12}
\partial_0^t : C_0^t (\mc B (G);V)^{U_{x_0}^{(r)}} \to V^{U_{x_0}^{(r)}} . 
\end{equation}
To see that the higher boundary maps are also continuous, we fix an $n+1$-polysimplex $\tau$
with faces $\tau_j$. The bilinear map of Fr\'echet spaces
\[
\begin{array}{ccc}
e_{U_{x_0}^{(r)}} \mc S (G) e_\tau \times \ep_\tau \otimes e_\tau V & \to &
\bigoplus_j e_{U_{x_0}^{(r)}} \mc S (G) e_{\tau_j} \widehat \otimes \ep_{\tau_j} \otimes
e_{\tau_j} V , \\ 
(f,\tau \otimes v) & \mapsto & f \widehat \otimes \partial (\tau) \otimes v =
\sum_j f \widehat \otimes [\tau : \tau_j] \tau_j \otimes v 
\end{array}
\]
is continuous, for it is made from the identity map of 
$e_{U_{x_0}^{(r)}} \mc S (G) e_\tau \widehat \otimes e_\tau V$, a finite linear combination
and the embeddings
\[
e_{U_{x_0}^{(r)}} \mc S (G) e_\tau \widehat \otimes e_\tau V \to
e_{U_{x_0}^{(r)}} \mc S (G) e_{\tau_j} \widehat \otimes e_{\tau_j} V .
\]
For $\sigma_i$ in the $G$-orbit of $\tau_j$
\[
e_{U_{x_0}^{(r)}} \mc S (G) e_{\tau_j} \widehat \otimes e_{\tau_j} V \to e_{U_{x_0}^{(r)}} 
\mc S (G) e_{\sigma_i} \underset{\mc H (G_{\sigma_i} ; U_{\sigma_i}^{(e)})}{\widehat \otimes} 
\ep_{\sigma_i} \otimes e_{\sigma_i} V 
\]
is a quotient map, so in particular continuous. Hence the composition
\[
e_{U_{x_0}^{(r)}} \mc S (G) e_\tau \times \ep_\tau \otimes e_\tau V \to 
C_n^t (\mc B (G);V)^{U_{x_0}^{(r)}}
\]
is also continuous. By construction the latter map is $\mc H (G_\tau,U_\tau^{(e)})$-balanced
and it extends $\partial_{n+1}$ on $e_{U_{x_0}^{(r)}} \mc H (G) e_\tau \times \ep_\tau 
\otimes e_\tau V$. Now the universal property of 
$\underset{\mc H (G_\tau ; U_\tau^{(e)})}{\widehat \otimes} $ says that there exists a 
unique continuous map
\begin{equation}\label{eq:3.13}
\partial_{n+1}^t : C_{n+1}^t (\mc B (G);V)^{U_{x_0}^{(r)}} \to 
C_n^t (\mc B (G);V)^{U_{x_0}^{(r)}}
\end{equation}
which extends $\partial_{n+1}$. 

By definition \eqref{eq:3.12} and \eqref{eq:3.13} are $\mc S (G,U_{x_0}^{(e)})$-linear.
As $C_n^t (\mc B (G);V)$ is endowed with the inductive limit topology, we can take the union 
over $r \in \N$ in \eqref{eq:3.12} and \eqref{eq:3.13} to find the required continuous maps 
$\partial_n^t$. These are homomorphisms of modules over 
$\bigcup_{r \in \N} \mc S (G,U_{x_0}^{(e)}) = \mc S (G). \qquad \Box$
\\[2mm]

Whereas Lemma \ref{lem:3.5} says that $\big( C_*^t (\mc B (G);V)^K ,\partial_*^t \big)$ is a 
differential complex in $\Mod_{Fr}(\mc S (G,K))$, the next lemma implies that it consists 
of projective objects.

\begin{lem}\label{lem:3.2}
Let $V \in \Mod_{Fr}(\mc S (G),K)$.
\enuma{
\item $C_n^t (\mc B (G);\mc S (G) e_K)^K$ is projective in 
$\Mod_{Fr}(\mc S(G,K) \widehat \otimes \mc S (G,K)^{op})$.
\item $C_n^t (\mc B (G);V)^K$ is projective in $\Mod_{Fr}(\mc S(G,K))$.
}
\end{lem}
\emph{Proof.} There are some technical complications, caused by the fact that \\
$C_n^t (\mc B (G);V)$ is not necessarily generated by its $K$-invariant vectors.

(a) Since $G$ acts transitively on the chambers of $\mc B (G)$, we may assume that $x_0$ and all the
$\sigma_i$ lie in the closure of a fixed chamber $c$. Then
\[
U_{x_0}^{(e+1)} \subset U_c^{(e+1)} \subset U_{\sigma_i}^{(e)} \; \text{for all } i=1,\ldots, d ,
\]
so abbreviating $U = U_{x_0}^{(e+1)}$ and $\sigma = \sigma_i$ we have
\[
e_U e_\sigma = e_\sigma = e_\sigma e_U .
\]
It follows that a typical direct summand of $C_n^t (\mc B (G);\mc S (G) e_U)^U$ looks like
\begin{multline*}
\big( \mc S(G,U) \widehat \otimes \mc S (G,U)^{op} \big) \sum_\rho 
(p_{\ep_\sigma\otimes\rho} \otimes p_{ \rho}) = S(G,U) e_\sigma \underset{\mc H (G_\sigma ; 
U_\sigma^{(e)})}{\widehat \otimes} \ep_\sigma \otimes e_\sigma \mc S(G,U) = \\ 
e_U S(G) e_\sigma \underset{\mc H (G_\sigma ; U_\sigma^{(e)})}{\widehat \otimes} 
\ep_\sigma \otimes e_\sigma \mc S(G) e_U .
\end{multline*}
These modules are projective in $\Mod_{Fr}(\mc S(G,U) \widehat \otimes \mc S (G,U)^{op})$,
because the left hand side is a direct summand of a free module. It follows readily from 
the definition of nuclearity that 
\begin{equation}\label{eq:3.7}
\mc S(G,U) \widehat \otimes \mc S (G,U)^{op} \cong \mc S( G \times G^{op} , U \times U^{op})
\end{equation}
as Fr\'echet algebras, see the proof of \cite[Lemma 1]{Mey2}.
By Proposition \ref{prop:3.1}, applied to $G \times G^{op}$, the functors
\begin{equation}\label{eq:3.4}
\begin{array}{ccc}
\Mod(\mc S (G \times G^{op},U \times U^{op})) & \longleftrightarrow & 
\Mod(\mc S (G \times G^{op},K \times K^{op})) \\
Y & \mapsto & Y^{K \times K^{op}} \\
e_{U \times U^{op}} \mc S (G \times G^{op}) e_{K \times K^{op}} Z & \leftarrow & Z
\end{array}
\end{equation}
provide an equivalence between Mod$(\mc S (G \times G^{op},K \times K^{op}))$ and a 
direct factor of Mod$(\mc S (G \times G^{op},U \times U^{op}))$. Because
$\mc S (G \times G^{op},U \times U^{op})$ is of finite rank as a module over
$\mc S (G \times G^{op},K \times K^{op})$, 
\[
e_{U \times U^{op}} \mc S (G \times G^{op}) e_{K \times K^{op}} Z = 
\mc S (G \times G^{op},U \times U^{op}) 
\underset{\mc S (G \times G^{op},K \times K^{op})}{\widehat \otimes} Z .
\]
It follows that the functors \eqref{eq:3.4} preserve the property Fr\'echet and preserve 
continuity of morphisms, so they remain equivalences for the appropriate categories of 
Fr\'echet modules. In particular $Y \mapsto Y^{K \times K^{op}}$ preserves projectivity, so
\[
\big( e_U S(G) e_\sigma \underset{\mc H (G_\sigma ; U_\sigma^{(e)})}{\widehat \otimes} 
\ep_\sigma \otimes e_\sigma \mc S(G) e_U \big)^{K \times K^{op}} =
e_K S(G) e_\sigma \underset{\mc H (G_\sigma ; U_\sigma^{(e)})}{\widehat \otimes} 
\ep_\sigma \otimes e_\sigma \mc S(G) e_K
\]
is projective in $\Mod_{Fr}(\mc S(G,K) \widehat \otimes \mc S (G,K)^{op})$. 
By \eqref{eq:3.2} $C_n^t (\mc B (G);\mc S (G) e_K)^K$ is a finite direct sum of
such modules.

(b) Apply the functor $\underset{\mc S (G,K)}{\widehat \otimes} V^K$ to part (a) and
use $V = \mc S (G) V^K . \qquad \Box$ 
\vspace{2mm}

\section{Continuity of the contraction}

We want to show that the contraction from Theorem \ref{thm:2.2} extends to 
$C_n^t (\mc B (G);V)^K$ by continuity. To that end we need a more concrete description
of $C_n^t (\mc B (G);V)^K$, at least in the universal case $V = \mc S (G)e_K$.

Let $\{\sigma_i\}$ a set of representatives for the $G$-orbits of polysimplices of 
$\mc B (G)$. We may assume that $x_0$ is among them and that all the $\sigma_i$ lie 
in a single chamber. We normalize the $G$-invariant metric $d$ on $\mc B (G)$ so that 
the diameter of a chamber is 1. Let $\ell : G \to \R_{\geq 0}$ be the length function
\[
\ell (g) = d (g x_0,x_0) + 1.
\]
As was shown in \cite[Section 9]{Vig1}, the topology on $\mc S (G)$ is defined by the 
norms
\[
p_m (f) = \norm{\ell^m f}_2 \qquad m \in \N.
\]
More precisely, $e_\tau \mc S (G) e_\sigma$ is the completion of 
$e_\tau \mc H (G) e_\sigma$ with respect to this collection of norms.

Note that the identification \eqref{eq:1.2} of the two appearances of 
$C_n (\mc B (G);\mc S(G)e_K)$ goes via the map
\begin{equation}\label{eq:3.6}
\begin{array}{ccc}
\alpha : \bigoplus\limits_{\sigma \in \mc B (G)^{(n)}} \C \sigma \otimes_\C
e_\sigma \mc S (G) e_K & \to & \bigoplus\limits_{i=1}^d \mc H (G) e_{\sigma_i}
\underset{\mc H (G_{\sigma_i})}{\otimes} \ep_{\sigma_i} \otimes e_{\sigma_i} \mc S (G) e_K \\
\sum_\sigma \sigma \otimes f_\sigma & \mapsto & \sum_\sigma e_\sigma g_\sigma^{-1} 
\otimes g_\sigma f_\sigma .
\end{array}
\end{equation}
Here we have chosen for each $\sigma$ an element $g_\sigma \in G$ such that 
$g_\sigma\sigma = \pm \sigma_i$, where $\sigma_i$ is the chosen representative of the 
$G$-orbit of $\sigma$. We fix such a choice of such $g_\sigma$ once and for all. 
Hence $g_\sigma f_\sigma \in e_{\sigma_i} \mc S (G) e_K$ for all $\sigma$. 
The fact that we tensor over $\mc H (G_{\sigma_i})$ makes the map $\alpha$ independent 
of the choices of the $g_\sigma$. The argument of \eqref{eq:2.3} shows that
\[
e_K \mc S (G) e_{\sigma_i} \underset{\mc H (G_{\sigma_i})}{\widehat \otimes} \ep_{\sigma_i} 
\otimes e_{\sigma_i} \mc S (G) e_K \cong \big( e_K \mc S (G) e_\sigma \widehat \otimes
(e_K \mc S (G) e_\sigma)^{op} \big) \sum_\rho (e_\rho \otimes e_{\ep_{\sigma_i} \rho} ) .
\]
As a Fr\'echet space this is a direct summand of
\begin{equation}\label{eq:3.8}
e_K \mc S (G) e_{\sigma_i} \widehat \otimes e_{\sigma_i} \mc S (G) e_K ,
\end{equation}
so the topology on $C_n^t (\mc B (G);\mc S (G) e_K)^K$ can be described with any 
defining family of seminorms on \eqref{eq:3.8}.

A general element of $C_n (\mc B(G) ;\mc S (G) e_K)$ can be written as 
$x = \sum_\sigma \sigma \otimes  f_\sigma$ where $f_\sigma \in e_{\sigma} \mc S (G) e_K$.  
We define $f_{-\sigma}$ by $f_{-\sigma} = -f_\sigma$ for all $\sigma$ (recall that we 
have fixed an orientation for all simplexes of $\mc B(G)$). We will use this notational 
convention from now on. Then $x$ is $K$-invariant if and only if $f_{k \sigma} = kf_\sigma$.
We define a norm $q_m \; (m \in \N)$ on $C_n (\mc B(G) ;\mc S (G) e_K)^K$ by
\begin{equation}\label{eq:4.5}
q_m \big( \sum\nolimits_\sigma \sigma \otimes f_\sigma \big) = 
\big( \sum\nolimits_\sigma \| (1 + d(\sigma,x_0) + \ell)^m g_\sigma f_\sigma \|_2^2 \big)^{1/2} .
\end{equation}
We remark at this point that this family of seminorms does depend on the choices of the 
elements $g_\sigma$, but not up to equivalence.

\begin{lem}\label{lem:4.1}
The Fr\'echet space $C_n^t (\mc B(G) ;\mc S (G) e_K)^K$ is the completion of\\
$C_n (\mc B (G);\mc H (G) e_K )^K$ with respect to the family of norms $\{ q_m \mid m \in \N \}$.
\end{lem}
\emph{Proof.}
By \cite[Section 9]{Vig1}, applied to $G \times G$, one defining family 
of norms on \eqref{eq:3.8} is
\[
q'_m ( f ) = \| (\ell_1 + \ell_2 )^m f \|_2 = \big( \int_{G \times G} 
(\ell (g_1) + \ell (g_2) )^ m |f(g_1,g_2)|^ 2 d g_1 d g_2 \big)^{1/2} .
\]
We retract $q'_m$ to $C_n (\mc B (G);\mc H (G) e_K )^K$ via \eqref{eq:3.6}. Notice that for
$g \in U_\sigma^{(e)} g_\sigma^{-1}$ the difference between $\ell(g)$ and $d(\sigma,x_0) + 1$ 
is at most $d(g x_0,g \sigma_i) \leq 1$, and hence inessential when it comes to these norms. 
Consider $x = \sum_\sigma \sigma \otimes f_\sigma \in C_n (\mc B (G);\mc H (G) e_K )^K$ and 
\begin{equation} \label{eq:4.4}
\alpha (x) = \sum\nolimits_\sigma e_\sigma g_\sigma^{-1} \otimes g_\sigma  f_\sigma .
\end{equation}
Since the right hand side of \eqref{eq:3.6} is a finite direct sum over the polysimplices 
$\sigma_i$, it suffices to consider the case that $x$ is supported on the $G$-orbit of one 
such $\sigma_i$. Then $e_\sigma g_\sigma^{-1}$ has support in $\{ g \in G : g \sigma_i = 
\sigma \}$, so the different $e_\sigma g_\sigma^{-1}$ have disjoint supports. Hence the sum 
\eqref{eq:4.4} is orthogonal for the $L_2$-norm, and this remains true if we multiply it with
the function $(\ell_1 + \ell_2 )^m$. The $L_2$-norm of $e_\sigma g_\sigma^{-1}$
is vol$\big( U_{\sigma_i}^{(e)} \big)^{-1/2}$, which is independent of $\sigma$.

It follows that $q'_m (\alpha (x))$ equals the right hand side of \eqref{eq:4.5}, up to a
constant factor. Consequently the norms $q_m$ with $m \in \N$ define the topology 
of $C_n^t (\mc B(G) ;\mc S (G) e_K)^K$. Now the result follows from the obvious density 
of $C_n (\mc B (G);\mc H (G) e_K )^K . \qquad \Box$ 
\\[1mm]

\begin{thm}\label{thm:4.2}
\enuma{
\item The differential complex 
\[
\mc S (G,K) \leftarrow C_*^t (\mc B (G);\mc S (G) e_K)^K
\]
is a projective resolution in $\Mod_{Fr}(\mc S (G,K) \widehat \otimes \mc S (G,K)^{op})$.
It admits a continuous contraction which is right $\mc S (G,K)$-linear.
\item Let $V \in \Mod_{Fr}(\mc S (G),K)$. Then
\[
V^K \leftarrow C_*^t (\mc B (G);V)^K \quad \text{and} \quad
V \leftarrow \mc S (G) 
e_K C_*^t (\mc B (G);V)
\]
are projective resolutions in $\Mod_{Fr}(\mc S (G,K))$ and in $\Mod_{Fr}(\mc S (G),K)$.
}
\end{thm}
\emph{Proof.}
(a) The projectivity was already established in Lemma \ref{lem:3.2}.\\
Like in the proof of Lemma \ref{lem:4.1} it suffices to consider an element  
$x = \sum_\sigma \sigma \otimes  f_\sigma \in C_n (\mc B (G);\mc H (G) e_K )^K$ . 
Recall from the text just above \eqref{eq:4.5} 
that the $K$-invariance is equivalent with $f_{k \sigma} = kf_{\sigma}$. Then 
\begin{equation}
\gamma_n (x) = \sum_{\tau \in \mc B(G)^{(n+1)}} \tau \otimes  
\sum_{\sigma \in \mc B (G)^{(n)}} \gamma_{\sigma,\tau} f_\sigma .
\end{equation}
By \eqref{eq:2.13} 
\[
F_\tau := \sum\nolimits_{\sigma \in \mc B (G)^{(n)}} \gamma_{\sigma,\tau}f_\sigma
\]
is invariant under $U_{\tau}^{(e)}$. By the $K$-equivariance of $\gamma$  
the element $\gamma_n(x)$ is $K$-invariant, and by our convention $F_{-\tau} = - F_\tau$. 
In view of Lemma \ref{lem:4.1} we have
\[
q_m (\gamma_n (x))^2 = \sum\nolimits_{\tau \in \mc B(G)^{(n+1)}}
\| (1 + d(\tau,x_0) + \ell )^m g_\tau F_\tau \|^2 .
\]
The $K$-invariance implies that  $F_{k\tau}=k F_\tau$ for all $\tau$ and $k\in K$. 
From the definition of $g_\tau$ it is clear that $g_{k\tau}k\in G_{\tau_i} g_\tau$. 
For all $h\in G_{\tau_i}$ and $g\in G$ we have $\ell(hg)\leq \ell(g)+2$. 
Since $d(k\tau,x_0)=d(\tau,x_0)$ for all $k\in K$, we obtain
\[
\| (1 + d(k\tau,x_0) + \ell )^m g_{k\tau} F_{k\tau} \|^2 \leq 
\| (3 + d(\tau,x_0) + \ell )^m g_\tau F_\tau \|^2 .
\]
Hence $q_m(\gamma_n(x))^2$ is bounded by 
\begin{multline*}
\sum_{\tau \in K \backslash \mc B(G)^{(n+1)}} [K : K \cap G_\tau] \,
\| (3 + d(\tau,x_0) + \ell )^m g_\tau F_\tau \|^2 = \\
\sum_{\tau \in K \backslash \mc B(G)^{(n+1)}} [K : K \cap G_\tau] \,
\big \| (3 + d(\tau,x_0) + \ell )^m \sum_{\sigma \in \mc B (G)^{(n)}} \gamma_{\sigma,\tau}
g_\tau f_\sigma \big\|^2_2 .
\end{multline*}
By Corollary \ref{cor:2.3}.c 
\begin{equation}\label{eq:4.6}
\leq M_\gamma^2 \sum_{\tau \in K \backslash \mc B(G)^{(n+1)}} \hspace{-5mm} [K : K \cap G_\tau]
\sum_{\sigma \in \mc B (G)^{(n)} : \gamma_{\sigma,\tau} \neq 0} \hspace{-5mm}
\big \| (3 + d(\tau,x_0) + \ell )^m g_\tau g_\sigma^{-1}(g_\sigma f_\sigma) \big\|^2_2 .
\end{equation}
The length function $\ell$ satisfies
\begin{align*}
\ell ( g_\tau g_\sigma^{-1} g) & \leq \ell (g_\tau g_\sigma^{-1}) + \ell (g) =
d (g_\tau g_\sigma^{-1} x_0,x_0) + 1 + \ell (g) \\
& = d (g_\sigma^{-1} x_0, g_\tau^{-1} x_0) + 1 + \ell (g) \\
& \leq d (\sigma,\tau) + d(g_\sigma^{-1} \sigma_i, g_\sigma^{-1} x_0) +
d(g_\tau^{-1} \tau_i, g_\tau^{-1} x_0) + 1 + \ell (g) \\
& \leq d(\sigma,\tau) + 3 + \ell (g) .
\end{align*}
Corollary \ref{cor:2.3}.b says that $\tau$ lies in the hull of $\sigma \cup \{x_0\}$ 
when $\gamma_{\sigma,\tau} \neq 0$, in which case 
\[
d(\sigma,\tau) + d (\tau,x_0) \leq 2 d (\sigma,x_0) + 2 d(\delta,x_0) 
\]
by Lemma \ref{lem:2.8}. We combine these length estimates to
\begin{multline}
3 + d(\tau,x_0) + \ell ( g_\tau g_\sigma^{-1} g) \leq 
6 + 2 d (\sigma,x_0) + 2 d(\delta,x_0) + \ell (g) \\
\leq \frac{6 + 2 d(x_0,\delta) + 1}{1 + 1} (1 + d(\sigma,x_0) + \ell (g))
:= c \, (1 + d(\sigma,x_0) + \ell (g)) .
\end{multline}
Therefore we may continue the estimate \eqref{eq:4.6} with
\begin{equation} \label{eq:4.7}
\leq M_\gamma^2 \sum_{\tau \in K \backslash \mc B(G)^{(n+1)}} \hspace{-5mm}
[K : K \cap G_\tau] \sum_{\sigma \in \mc B (G)^{(n)} : \gamma_{\sigma,\tau} \neq 0} \hspace{-5mm}
\big \| c^m (1 + d(\sigma,x_0) + \ell )^m g_\sigma f_\sigma \big\|^2_2 .
\end{equation}
By Corollary \ref{cor:2.3}.b the sets $S_\tau:=\{ \sigma \in \mc B (G)^{(n)} : 
\gamma_{\sigma,\tau} \neq 0 \}$ and $kS_\tau=\{\sigma \in \mc B (G)^{(n)} : 
\gamma_{\sigma,k \tau} \neq 0 \}$ are disjoint if $k \tau \neq \pm \tau$.
By Lemma \ref{lem:2.8} their union is contained in 
\[
\mc B (G)^{(n)}_\tau := \{\sigma \in \mc B (G)^{(n)} : 
d(\sigma,x_0) \geq d(\tau,x_0) - d (\delta,x_0) \} ,
\]
so \eqref{eq:4.7} is bounded by (using the invariance $f_{k\sigma}=k f_\sigma$)
\begin{equation}\label{eq:4.8}
\leq M_\gamma^2 c^{2m} \sum_{\tau \in K \backslash \mc B(G)^{(n+1)}} 
\sum_{\sigma \in \mc B (G)^{(n)}_\tau} 
\big \| (1 + d(\sigma,x_0) + \ell )^m g_\sigma f_\sigma \big\|^2_2 .
\end{equation}
By the Cartan decomposition $G_{x_0} \backslash \mc B (G)$ is in bijection with a Weyl
chamber in an apartment $A$. As $K$ is of finite index in $G_{x_0}$, this shows that
$K \backslash \mc B(G)^{(d)}$ is of polynomial growth. Choose $N \in 2\N$ such that
\[
b_N := \sum\nolimits_{K \backslash \mc B(G)^{(n+1)}} (1 + d(\tau,x_0) )^{-N} \text{ is finite.} 
\]
This enables us to estimate \eqref{eq:4.8} by 
\begin{align*}
\leq M_\gamma^2 c^{2m} \sum_{\tau \in K \backslash \mc B(G)^{(n+1)}} 
(1 + d(\tau,x_0) )^{-N} \sum_{\sigma \in \mc B (G)^{(n)}_\tau}
\big \| (1 + d(\sigma,x_0) + \ell )^{m+N/2} g_\sigma f_\sigma \big\|^2_2 \\ 
\leq M_\gamma^2 c^{2m} b_N \sum_{\sigma \in \mc B (G)^{(n)}} 
\big \| (1 + d(\sigma,x_0) + \ell )^{m+N/2} g_\sigma f_\sigma \big\|^2_2 \; = \;
M_\gamma^2 c^{2m} b_N q_{m+N/2} (x)^2 .
\end{align*}
Altogether we obtained
\begin{equation}
q_m (\gamma_n (x)) \leq M_\gamma c^m \sqrt{b_N} q_{m+N/2}(x) ,
\end{equation}
from which we conclude that on $C_n (\mc B (G);\mc S (G) e_K)^K$ the map $\gamma_n$ is continuous
with respect to the family of norms $\{ q_m \mid m \in \N\}$. From Theorem \ref{thm:2.2} we
know that $\gamma_n$ is right $\mc H (G,K)$-linear, so by Lemma \ref{lem:4.1} it extends 
continuously to a right-$\mc S (G,K)$-linear map 
\[
\gamma_n^t: C_n^t (\mc B (G);\mc S (G) e_K)^K \to C_{n+1}^t (\mc B (G);\mc S (G) e_K)^K .
\]
The relation $\partial_{n+1} \gamma_n + \gamma_{n-1} \partial_n =$ id extends by continuity
to $\partial_{n+1}^t \gamma_n^t + \gamma_{n-1}^t \partial_n^t =$ id.\\
(b) The first statement follows from (a) upon applying the functor 
$\underset{\mc S (G,K)}{\widehat \otimes} V^K$. The second is a consequence of the first
and Corollary \ref{cor:3.6}. $\qquad \Box$ 
\\[2mm]

We remark that Theorem \ref{thm:4.2} does not imply that $V \leftarrow C_*^t (\mc B (G);V)$
is a resolution. Although this is true, one needs more sophisticated techniques to
prove it -- see the next section.
The main use of Theorem \ref{thm:4.2} is to compute and compare Ext-groups:

\begin{prop}\label{prop:4.3}
Let $V,W \in \Mod(\mc S (G))$ with $V$ admissible. 
\enuma{
\item There is a natural isomorphism
\[
\mr{Ext}^n_{\mc H (G)}(V,W) \cong \mr{Ext}^n_{\mc S (G)}(V,W) .
\]
If $W \in \Mod_{LF}(\mc S (G))$, then these are also isomorphic to
$\mr{Ext}^n_{\Mod_{LF} (\mc S (G))} (V,W)$.
\item Suppose that moreover $V,W \in \Mod(\mc S (G),K)$. There are natural isomorphisms
\begin{multline*} 
\mr{Ext}^n_{\mc H (G)} (V,W) \cong \mr{Ext}^n_{\mc H (G,K)} (V^K,W^K) \\
\cong \mr{Ext}^n_{\mc S (G,K)} (V^K,W^K) \cong \mr{Ext}^n_{\Mod(\mc S (G),K)} (V,W) .
\end{multline*}
\item If furthermore $V,W \in \Mod_{Fr}(\mc S (G),K)$, 
then the groups from (b) are also naturally isomorphic to
\[
\mr{Ext}^n_{\Mod_{Fr} (\mc S (G,K))} (V^K,W^K) \quad \text{and} \quad
\mr{Ext}^n_{\Mod_{Fr} (\mc S (G),K)} (V,W) .
\]
}
\end{prop}
\emph{Proof.} (b)
The outer isomorphisms follow from Theorem \ref{thm:2.1} and 
Corollary \ref{cor:3.6}. For the middle one, we observe that by Theorem \ref{thm:2.2}.a
\begin{align}
\nonumber \mr{Ext}^n_{\mc H (G,K)} (V^K & ,W^K) = H^n \Big( \mr{Hom}_{\mc H (G,K)} \big( 
C_* (\mc B (G);V)^K, W^K \big), \mr{Hom}(\partial_*,W^K ) \Big) \\
\label{eq:4.2} = H^n & \Big( \mr{Hom}_{\mc S (G,K)} \big( \mc S (G,K) 
\underset{\mc H (G,K)}{\otimes} C_* (\mc B (G);V)^K, W^K \big), \mr{Hom}(\partial_*,W^K) \Big) .
\end{align}
As $V$ is admissible, 
\[
\mc S (G,K) \underset{\mc H (G,K)}{\otimes} C_* (\mc B (G);V)^K = C_*^t (\mc B (G);V)^K .
\] 
By Theorem \ref{thm:1.2}, Frobenius reciprocity and Lemma \ref{lem:3.2}.b this module
is projective in Mod$(\mc S (G,K)$ and in $\Mod_{Fr}(\mc S (G,K))$. Moreover it is finitely 
generated, so every module map to a Fr\'echet $\mc S (G,K)$-module is automatically continuous. 
Therefore \eqref{eq:4.2} equals
\[
H^n \Big( \mr{Hom}_{\mc S (G,K)} \big( C_*^t (\mc B (G);V)^K, W^K \big), 
\mr{Hom}(\partial_*,W^K) \Big) ,
\]
which by Theorem \ref{thm:4.2} is $\mr{Ext}^n_{\mc S (G,K)} (V^K,W^K)$.\\ 
(c) In case $W^K \in \Mod_{Fr}(\mc S (G,K))$, the above argument also shows that we obtain 
the same answer if we work in $\Mod_{Fr}(\mc S (G,K))$. By Corollary \ref{cor:3.6}, these
Ext-groups are naturally isomorphic to $\mr{Ext}^n_{\Mod_{Fr} (\mc S (G),K)} (V,W)$. \\ 
(a) The first statement was proven in \cite[Section 9]{ScZi2}, using the results of
Meyer \cite{Mey2}. Here we provide an alternative proof. Recall that the Bernstein
decomposition of Mod$(\mc H (G))$ is given by idempotents in the centre of the
category Mod$(\mc H (G))$ \cite{BeDe}. Hence $\Mod(\mc S (G))$ and 
$\Mod_{LF}(\mc S (G))$ admit an analogous decomposition. This persists to
Ext-groups, so we may and will assume that $V$ and $W$ live in a single Bernstein
component $\Omega$. Choose $e \in \N$ such that all representations in $\Omega$
are generated by their $U_{x_0}^{(e)}$-invariant vectors. Then $V,W \in \Mod
(\mc S (G), U_{x_0}^{(e)})$ and moreover $W^K \in \Mod_{Fr}(\mc S (G,K))$
if $W \in \Mod_{LF}(\mc S (G))$. Now we can apply parts (b) and (c).
$\qquad \Box$ \\[2mm]

The admissibility of $V$ is necessary in Proposition \ref{prop:4.3}.
The difference can already be observed in degree $n = 0 \!: \mr{Ext}^0_{\mc S (G)}(V,W)
= \mr{Hom}_{\mc S (G)} (V,W)$ can be smaller than $\mr{Ext}^0_{\mc H (G)}(V,W)
= \mr{Hom}_{\mc H (G)} (V,W)$ for general $V,W \in \Mod(\mc S (G))$. In 
$\Mod_{LF}(\mc S (G))$ we usually get an even smaller space of morphisms, because they 
are required to be continuous.
\vspace{1mm}

\section{Bornological modules}

The content of Sections 3 and 4 can be formulated nicely with bornologies. In this
section we work in the category $\Mod_{bor}(A)$ of complete bornological
modules over a bornological algebra $A$, as in \cite{Mey1}. The corresponding
tensor product is the completed bornological tensor product over $A$, which we denote
by $\widetilde \otimes_A$. In case $A = \C$, we suppress it from the notation.

We endow $\mc H (G)$ with the fine bornology, so a subset of $\mc H (G)$ is considered
to be bounded if it is contained in a finite dimensional linear subspace of $\mc H (G)$
and over there is bounded in the usual sense. On $\mc S (G)$ we use the precompact
bornology, which means that a subset is bounded if and only if it is contained in
some compact subset. Since we use the inductive limit topology on $\mc S (G)$, every
bounded set is contained in a compact subset of $\mc S (G,U)$ for some compact open
subgroup $U \subset G$.

By \cite[Lemma 2]{Mey2} we have
\[
\mc H (G) \widetilde \otimes \mc H (G) \cong \mc H (G \times G) \quad \text{and}
\quad  \mc S (G) \widetilde \otimes \mc S (G) \cong \mc S (G \times G).
\]
These isomorphisms, the second of which does not hold for the algebraic or the completed
projective tensor product, to some extent explain why bornology is a convenient technique
in our situation.

Since $\mc H (G)$ has the fine bornology, the projectivity properties from Theorems
\ref{thm:1.2} and \ref{thm:2.2} carry over. Hence for any 
$V \in \Mod_{bor}(\mc H (G))$ the modules
\begin{equation}
\begin{array}{lll}
C_n (\mc B (G);V) & \in & \Mod_{bor}(\mc H (G)) , \\
C_n (\mc B (G);\mc H (G)) & \in & 
\Mod_{bor}(\mc H (G) \widetilde \otimes \mc H (G)^{op}) , \\
C_n (\mc B (G);V)^K & \in & \Mod_{bor}(\mc H (G,K)) , \\
C_n (\mc B (G);\mc H (G) e_K)^K & \in & 
\Mod_{bor}(\mc H (G,K) \widetilde \otimes \mc H (G,K)^{op})
\end{array} 
\end{equation}
are projective in the respective categories.

The categories of topological $\mc S (G)$-modules that we used in the previous sections
are full subcategories of $\Mod_{bor}(\mc S (G))$. To see this, we endow all $V,W \in
\Mod_{LF}(\mc S (G))$ with the precompact bornology. Any $\mc S (G)$-module map
$\phi : V \to W$ sends $V^U$ to $W^U$, for any compact open subgroup $U$. By the
definition of the inductive limit topology, $\phi$ is continuous if and only if
$\phi \big|_{V^U}$ is continuous for all $U$. Since $V^U$ and $W^U$ are Fr\'echet spaces,
the latter condition is equivalent to boundedness of $\phi \big|_{V^U}$. As the bornology
on $V$ is the inductive limit of the bornologies on the $V^U$, this in turn is equivalent
to boundedness of $\phi$. Hence 
\[
\mr{Hom}_{\Mod_{bor} (\mc S (G))}(V,W) = \mr{Hom}_{\Mod_{LF} (\mc S (G))}(V,W) .
\]
Since bornological tensor products commute with inductive limits, the definition 
\eqref{eq:3.2} can be simplified to 
\[
C_n^t (\mc B (G);V) = \bigoplus_{i=1}^d \mc S (G) \underset{\mc H (G_{\sigma_i} ;
U_{\sigma_i}^{(e)})}{\widetilde \otimes} \ep_{\sigma_i} \otimes V^{U_{\sigma_i}^{(e)}}
\quad \text{for} \quad V \in \Mod_{bor}(\mc S (G)) .
\]
The same argument as in the proof of Theorem \ref{thm:1.2} shows that this is a 
projective object of $\Mod_{bor}(\mc S (G))$. By \eqref{eq:2.3}
\[
C_n^t (\mc B (G);\mc S (G)) = \bigoplus_{i=1}^d \big( \mc S (G) \widetilde \otimes
\mc S (G)^{op} \big) \sum_\rho (p_\rho \otimes p_{\ep_{\sigma_i} \rho})
\]
is projective in $\Mod_{bor}(\mc S (G) \widetilde \otimes \mc S (G)^{op}) = \Mod_{bor}
(\mc S (G \times G^{op}))$. Just like for Fr\'echet modules, Proposition \ref{prop:3.1} 
remains valid for bornological modules. It follows that 
\begin{align}\label{eq:5.2}
& C_n^t (\mc B (G);V)^K = \bigoplus_{i=1}^d e_K \mc S (G) \underset{\mc H (G_{\sigma_i} ;
U_{\sigma_i}^{(e)})}{\widetilde \otimes} \ep_{\sigma_i} \otimes e_{U_{\sigma_i}^{(e)}} V
\in \Mod_{bor}(\mc S (G,K)) , \\
& C_n^t (\mc B (G);\mc S (G) e_K )^K \in 
\Mod_{bor}(\mc S (G,K) \widetilde \otimes \mc S (G,K)^{op}) \nonumber
\end{align}
are projective. Furthermore we note that, by the associativity of 
bornological tensor products and by \eqref{eq:3.1},
\begin{align}\label{eq:5.1} 
C_n^t (\mc B (G);V) \; & = \; \mc S (G) \underset{\mc H (G)}{\widetilde \otimes} 
C_n (\mc B (G);V) , \\ 
\nonumber C_n^t (\mc B (G);V)^K \; = \; e_K \mc S (G) 
\underset{\mc H (G)}{\widetilde \otimes} C_n (\mc B (G);V)^K \; & = \; \mc S (G,K) 
\underset{\mc H (G,K)}{\widetilde \otimes} C_n (\mc B (G);V)^K . 
\end{align}
According to \cite[Theorem 22]{Mey2} the embedding of bornological algebras $\mc H (G) 
\to \mc S (G)$ is isocohomological. Together with Theorem \ref{thm:1.2}.b and 
\cite[Theorem 35]{Mey1} this implies that, for any $V \in \Mod_{bor} (\mc S (G))$ which is
generated by its $U_{x_0}^{(e)}$-invariant vectors for some $e \in \N$,
\begin{equation} \label{eq:5.3}
\big( \mc S (G) \underset{\mc H (G)}{\widetilde \otimes} C_* (\mc B (G);V), \partial_*^t \big) 
= \big( C_*^t (\mc B (G);V) , \partial_*^t \big)
\end{equation}
is a resolution of $V$ in $\Mod_{bor} (\mc S (G))$.

\begin{thm}\label{thm:5.1}
Let $V,W \in \Mod_{bor}(\mc S (G))$.
\enuma{
\item There is a natural isomorphism
\[
\mr{Ext}^n_{\Mod_{bor}(\mc H (G))}(V,W) \cong
\mr{Ext}^n_{\Mod_{bor}(\mc S (G))}(V,W) . 
\]
\item Suppose that moreover $V,W \in \Mod (G,K)$. There are natural isomorphisms
\begin{multline*} 
\mr{Ext}^n_{\Mod_{bor}(\mc H (G))} (V,W) \cong 
\mr{Ext}^n_{\Mod_{bor}(\mc H (G,K))} (V^K,W^K) \\
\cong \mr{Ext}^n_{\Mod_{bor}(\mc S (G,K))} (V^K,W^K) 
\cong \mr{Ext}^n_{\Mod_{bor}(\mc S (G))} (V,W) .
\end{multline*}
} 
\end{thm}
\emph{Proof.} 
Of course this is a straightforward consequence of Meyer's result \eqref{eq:5.3}. 
Even Theorem \ref{thm:2.2} is not really needed, only the existence of some projective
resolution. Here we show how the theorem can be derived from Theorem \ref{thm:4.2}. \\
(b) The outer isomorphisms follow from the bornological versions of Theorem \ref{thm:2.1}
and Proposition \ref{prop:3.1}. As concerns the middle one, by Theorem \ref{thm:2.2}.a
\begin{multline*}
\mr{Ext}^n_{\Mod_{bor}(\mc H (G,K))} (V^K,W^K) = \\
H^n \big( \mr{Hom}_{\Mod_{bor}(\mc H (G,K))} \big( C_* (\mc B (G);V)^K , W^K \big) ,
\mr{Hom}_{bor}(\partial_* , W^K) \big) .
\end{multline*}
By \eqref{eq:5.1} and by Frobenius reciprocity this is isomorphic to
\[
H^n \big( \mr{Hom}_{\Mod_{bor}(\mc S (G,K))} \big( C_*^t (\mc B (G);V)^K , W^K \big) ,
\mr{Hom}_{bor}(\partial_*^t , W^K) \big) ,
\]
which by Theorem \ref{thm:4.2} and \eqref{eq:5.2} is 
$\mr{Ext}^n_{\Mod_{bor}(\mc S (G,K))} (V^K,W^K)$.\\
(a) By the same argument as for Proposition \ref{prop:4.3}, this follows from part (b).
$\qquad \Box$ \\[2mm]

We remark that in general 
\[
\mr{Ext}^n_{\Mod_{bor}(\mc H (G))} (V,W) \not\cong \mr{Ext}^n_{\mc H (G)} (V,W) .
\]
The reason is that morphisms in $\Mod_{bor}(\mc H (G))$ have to be bounded, which is
a nontrivial condition if $V$ is not admissible.
\vspace{2mm}

\section{Generalization to disconnected reductive groups}

In the final section we take a more general point of view, we let $G = \mc G (\mh F)$
be an algebraic group whose identity component $G^\circ = \mc G^\circ (\mh F)$ is
linear and reductive. We will show how the results of the previous sections can be
generalized to such groups.

First we discuss the categorical issues. Since $G$ acts on $G^\circ$ by conjugation,
it acts on $\Mod (G^\circ) = \Mod (\mc H (G^\circ ;\C))$ and on the centre of this 
category. If $e_\Omega^\circ$ is a central idempotent of $\Mod (G^\circ)$, then
\[
e_\Omega := \sum_{g \in G / \mr{Stab}_G (\Omega^\circ)} g e_\Omega^\circ g^{-1} 
\]
is a central idempotent of Mod$(G)$. It follows that the category of smooth 
$G$-representations on complex vector spaces admits a factorization, parametrized by the 
$G$-orbits of Bernstein components of $G^\circ$:
\begin{equation}\label{eq:6.2}
\Mod (G) = \prod_{\Omega = G \Omega^\circ / G} \Mod_\Omega (G) = 
\prod_{\Omega = G \Omega^\circ / G} e_\Omega \Mod (G) .
\end{equation}
However, in contrast with the Bernstein decomposition for connected reductive
$p$-adic groups, it is possible that $\Mod_\Omega (G)$ is decomposable.

Following \cite{BuKu} we call a compact open subgroup $U \subset G$ (or more
precisely the idempotent $e_U \in \mc H (G)$) a type if the category $\Mod (G,U)$
is closed under the formation of subquotients in $\Mod (G)$. In case $\mc G$ is 
connected, \cite[Proposition 3.6]{BuKu} shows that these are precisely the compact
open subgroups for which Theorem \ref{thm:2.1} holds.

\begin{lem}\label{lem:6.1}
Let $U \subset G^\circ$ be a type for $G^\circ$. Then Theorem \ref{thm:2.1}
holds for $(G,U)$. In particular $V \mapsto V^U$ defines an equivalence of
categories $\Mod (G,U) \to \Mod (\mc H (G,U))$ and $\Mod (G,U)$ is a direct
factor of $\Mod (G)$. 
\end{lem}
\emph{Proof.}
Part (a) of Theorem \ref{thm:2.1} is \cite[Proposition 3.3]{BuKu}. 
By \cite[Proposition 3.6]{BuKu} there exists a finite collection $\Lambda$ of
Bernstein components for $G^\circ$ such that
\[
\Mod (G^\circ, U) = \prod_{\Omega^\circ \in \Lambda} \Mod_{\Omega^\circ} (G^\circ) . 
\]
We claim that
\begin{equation}\label{eq:6.3}
\Mod (G, U) = \prod_{\Omega \in G \Lambda / G} \Mod_{\Omega} (G)  .
\end{equation}
First we consider $V \in \Mod (G,U)$ as a $\mc H (G^\circ)$-module. As such
\[
V = \mc H (G) V^U = \sum_{g \in G / G^\circ} \mc H (G^\circ) g V^U =
\sum_{g \in G / G^\circ} \mc H (G^\circ) V^{g U g^{-1}} ,
\]
so $V$ is generated by $\sum_{g \in G / G^\circ} V^{g U g^{-1}}$. Since
\[
\mc H (G^\circ) V^{g U g^{-1}} \in 
\prod_{\Omega^\circ \in g \Lambda} \Mod_{\Omega^\circ} (G^\circ) = 
\prod_{\Omega^\circ \in g \Lambda} e_{\Omega^\circ} \Mod (G^\circ) ,
\]
$V$ lies in 
\begin{multline*}
\prod_{\Omega^\circ \in G \Lambda} e_{\Omega^\circ} \Mod (G^\circ) =
\prod_{\Omega^\circ \in G \Lambda / G} \sum_{g \in G / \mr{Stab}_G (\Omega^\circ)} 
g e_{\Omega^\circ} g^{-1} \Mod (G^\circ) \\
= \prod_{\Omega = G \Omega^\circ / G \in G \Lambda / G} e_\Omega \Mod (G^\circ) .
\end{multline*}
In particular $V \in \prod_{\Omega \in G \Lambda / G} \Mod_\Omega (G)$.
Conversely, let $W \in \Mod_\Omega (G)$ with $\Omega = G \Omega^\circ / G \in 
G \Lambda /G$. Then
\[
W = e_\Omega W = \sum_{g \in G / \mr{Stab}_G (\Omega^\circ)} g e_\Omega^\circ g^{-1} W
\in \prod_{G / \mr{Stab}_G (\Omega^\circ)} \Mod_{g \Omega^\circ}(G^\circ) .
\]
As $\Omega^\circ \in \Lambda$, the latter category is a direct factor of
$\prod_{g \in G / N_G (U)} \Mod (G^\circ, g U g^{-1})$. Now we can write
\[
W = \sum_{g \in G / N_G (U)} \mc H (G^\circ) W^{g U g^{-1}} = 
\sum_{g \in G / N_G (U)} \mc H (G^\circ) g W^U = \mc H (G) W^U ,
\]
which verifies our claim \eqref{eq:6.3}. Finally, by \eqref{eq:6.2} 
$\prod_{\Omega \in G \Lambda / G} \Mod_\Omega (G)$ is a direct factor of 
$\Mod (G). \qquad \Box$ 
\\[2mm]

Because the decomposition \eqref{eq:6.2} is defined in terms of central idempotents
of $\Mod (G)$, it can be handled in the same way as the usual Bernstein
decomposition. In particular the proofs of \cite[Lemmas 2.2 and 2.3]{ScZi1} remain
valid. Using these proofs in various categories leads to:

\begin{cor}\label{cor:6.3}
Let $U \subset G^\circ$ be a type for $G^\circ$. Then Theorem \ref{thm:2.1} holds
for $(G,U)$ in $\Mod (\mc S (G)), \Mod_{LF}(\mc S (G)), \Mod_{bor}(\mc S (G))$ and
$\Mod_{bor}(\mc H (G))$.
\end{cor}

Next we deal with the affine building of $G$. As a metric space, it is defined as
\[
\mc B (G) = \mc B (\mc G ,\mh F) := \mc B (\mc G^\circ ,\mh F) =
\mc B (G^\circ / Z(G^\circ)) \times X_* (Z (G^\circ)) \otimes_\Z \R .
\]
The action of $G^\circ$ on $\mc B (G^\circ / Z (G^\circ))$ is extended to $G$ in
the following way. There is a bijection between maximal compact subgroups of
$G^\circ / Z (G^\circ)$ and vertices of $\mc B (G^\circ / Z (G^\circ))$, which
associates to a vertex $x$ its stabilizer $K_x$. For any $g \in G$ the subgroup
$g K_x g^{-1} \subset G^\circ / Z(G^\circ)$ is again maximal compact, so of the form
$K_y$ for a unique vertex $y$ of $\mc B (G^\circ / Z (G^\circ))$. We define 
$g(x) = y$ and extend this by interpolation to an isometry of $\mc B (G^\circ / 
Z (G^\circ))$.

Since $Z(G^\circ)$ is a characteristic subgroup of $G^\circ$ and $G^\circ$ is normal
in $G ,\; G$ acts on $Z(G^\circ)$ by conjugation. This induces an action of $G$ on
$X_* (Z(G^\circ))$ which extends the action of $G^\circ$. 

Because the polysimplicial structure on $\mc B (G^\circ / Z (G^\circ))$ is natural,
it is preserved by $G$. Unfortunately, no such thing holds for $X_* (G^\circ / 
Z(G^\circ)) \otimes_\Z \R$, so our choice of a polysimplicial structure is in general
not stable under the $G$-action. Even worse, if $G$ acts in a complicated way on
$X_* (G^\circ / Z(G^\circ))$, it can be very difficult to find a suitable 
polysimplicial structure on $X_* (G^\circ / Z(G^\circ)) \otimes_\Z \R$. A serious
investigation of this problem would lead us quite far away from the theme of the
paper, so we avoid it by means of the following assumption.

\begin{cond}\label{cond:6.4}
There exists a $G$-stable root system of full rank in $X_* (Z(G^\circ))$. 
\end{cond}

Under this condition the affine Coxeter complex of the root system is a suitable
polysimplicial structure on $X_* (G^\circ / Z(G^\circ)) \otimes_\Z \R$. In most
examples $G / G^\circ$ is small and the condition is easily seen to be 
fulfilled.

We also need a slightly improved version of the groups $U_\sigma^{(e)}$. To define
it, we have to go through a part of the construction from \cite{ScSt}. Let 
$T = \mc T (\mh F)$ be a maximal $\mh F$-split torus of $G^\circ$ and let $\Phi$
be the root system of $(G^\circ ,T)$. Furthermore denote by $U^+$ and $U^-$ the
unipotent subgroups of $G^\circ$ corresponding to some choice of positive and
negative roots. The new group $U_\sigma^{[e]}$ will admit a factorization like
\eqref{eq:2.4}. 

First we assume that $\mc G^\circ$ is quasi-split over $\mh F$, so that 
$Z_{G^\circ}(T)$ is a maximal torus of $G^\circ$. We can keep $U_\sigma^{(e)} 
\cap U^+$ and $U_\sigma^{(e)} \cap U^-$, but we have to change $U_\sigma^{(e)}
\cap Z_{G^\circ} (T)$. Let $Z_{G^\circ}(T)^{mc}_r \; (r \in \R_{\geq 0})$ be the
``minimal congruent filtration'' of the torus $Z_{G^\circ}(T)$, as defined in
\cite[\S 5]{Yu}, and put 
\[
Z_{G^\circ}(T)^{mc}_{e+} = \bigcup\nolimits_{r > e} Z_{G^\circ}(T)^{mc}_r .
\]
Following Yu we define, for $e \in \N$ and a polysimplex $\sigma$ of $\mc B(G) =
\mc B (\mc G^\circ ,\mh F)$:
\[
U_\sigma^{[e]} := (U_\sigma^{(e)} \cap U^+) \, Z_{G^\circ}(T)^{mc}_{e+} \,
(U_\sigma^{(e)} \cap U^-) .
\]
For general (not quasi-split) $\mc G^\circ$ the subgroups $U_\sigma^{[e]} \subset 
G^\circ$ are obtained from those in the quasi-split case via etale descent, as in 
\cite[Section 1.2]{ScSt}. For this system of subgroups Proposition \ref{prop:1.1}
and all the results of \cite{ScSt} hold. 

According to \cite[\S 9.4]{Yu}, there exists an affine group scheme 
$\mc G_\sigma^\circ$ such that $\mc G_\sigma^\circ (\mc O) = G_\sigma^\circ$ is 
the pointwise stabilizer of $\sigma$ and 
\[
U_\sigma^{[e]} = \ker \big(\mc G_\sigma^\circ (\mc O) \to 
\mc G_\sigma^\circ (\mc O / \pi^{e+1} \mc O) \big) .
\]
Consequently $U_\sigma^{[e]}$ is stable under any automorphism of the affine group 
scheme $\mc G_\sigma^\circ$, which is not guaranteed in full generality in 
\cite{ScSt}. Clearly this applies to the action of $N_G (G_\sigma^\circ)$ on
$\mc G_\sigma^\circ$, so $g U_\sigma^{[e]} g^{-1} = U_\sigma^{[e]}$ for all
$g \in N_G (G_\sigma^\circ)$. Proposition \ref{prop:1.1}.d and the definition of
the action of $G$ on $\mc B (G)$ then entail
\[
g U_\sigma^{[e]} g^{-1} = U_{g \sigma}^{[e]} \quad \text{for all } g \in G . 
\]
Therefore Proposition \ref{prop:1.1} holds for $G$ with the system of subgroups 
$U_\sigma^{[e]}$, which is enough to make everything from \cite{MeSo1} work.

With the above adjustments almost everything in the preceeding sections generalizes to
disconnected reductive groups, only the proof of Theorem \ref{thm:2.2}.a needs 
a little more care. It is here that we use the condition \ref{cond:6.4}. The
construction of the contraction $\gamma$ of $C_* (A;\Z)$ in \cite[Section 2.1]{OpSo1} 
applies to an apartment $A$ spanned by an integral root system. Thus the assumed 
root system in $X_* (Z(G^\circ))$, together with the roots of $(G^\circ,T)$, 
functions as a book-keeping device to write down a contraction which has
some nice properties. The remainder of the proof of Theorem \ref{thm:2.2}.a needs
no modification.

We conclude that:
\begin{thm}\label{thm:6.5}
Under Condition \ref{cond:6.4} all the results of Sections 1--5 are valid for
$G$ with the system of subgroups $U_\sigma^{[e]}$. 
\end{thm}

\end{document}